\documentclass{amsart}
\usepackage{amsmath}
\usepackage{amssymb}
\usepackage[dvips]{graphicx}
\usepackage[latin1]{inputenc}
\usepackage{amssymb,latexsym}

\newtheorem{theorem}{Theorem}
\newtheorem{corollary}{Corollary}

\newtheorem{proposition}{Proposition}

\title[Ricci curvature and conformality to
spheres ]{Ricci curvature and conformality of Riemannian  manifolds
to spheres }

\author[Najoua Gamara, Salem Eljazi and Habiba Guemri ]{}


\email{najoua.gamara@fst.rnu.tn,  ngamara7@gmail.com}
\email{eljazisalem@gmail.com} \email{habiba.guemri@gmail.com}

\begin{document}
\maketitle

\centerline{ Najoua Gamara, Salem Eljazi }
\medskip
{\footnotesize
  \centerline{Facult\'e des Sciences de
Tunis, D\'epartement de Math\'ematiques} \centerline{Campus universitaire -El Manar II- 2092, Tunis, Tunisie.}}


\medskip
\centerline{Habiba Guemri}
\medskip
{\footnotesize \centerline{Institut Préparatoire aux Etudes d'Ingénieurs d'El Mannar}
\centerline{Campus universitaire -El Manar II - 2092, Tunis, Tunisie.}
\smallskip
}
\medskip
\bigskip

\begin{abstract}
In this paper we give bounds on the least eigenvalue of the
conformal Laplacian and the Yamabe invariant of a compact Riemannian
manifold in terms of the Ricci curvature and the diameter
and deduce a sufficient condition for the manifold to be conformally equivalent to a
sphere.\\\\ R\'esum\'e: Soit $(M,g)$ une vari\'et\'e riemannienne
compacte sans bord de dimension $n$. En utilisant des bornes
inf\'erieures sur la courbure de Ricci et le diam\`etre de $(M,g)$, on minore la
plus petite valeur propre du laplacien conforme ainsi que
l'invariant de Yamabe de cette vari\'et\'e. On en d\'eduit certaines
conditions pour que $(M,g)$ soit conform\'ement diff\'eomorphe \`a
la sph\`ere unit\'e de m\^eme dimension.\\\\
\end{abstract}

{\bf Motivation }\\ \indent The purpose of the paper is to give
conditions on some topological or geometrical invariants of a smooth
compact Riemannian manifold $(M,g)$ without boundary, to be
conformally diffeomorphic to the sphere of the same dimension
equipped with its canonical metric. This problem has an interesting
history. Indeed, the first approach was based on the use of the
conformal automorphisms group of the manifold denoted by $C(M,g)$. It
was shown that the non-compactness of the connected component of the
identity in $C(M,g)$ implies that $(M,g)$ is conformally equivalent
to the sphere $S^{n}$ for $n\geq 3$ (this is due to M.Obata,
 \cite{Ob1} and  \cite{Ob2}). Unfortunately, there was a gap in
\cite{Ob2} involving Obata's use of a certain theorem and
K.R.Gutschera in \cite{Gu} gave some counterexamples and finally
J.Lafontaine \cite{La} completed the proof in 1988. The compactness
of the whole group $C(M,g)$ was shown by J. Ferrand \cite{Fe}. An
alternate approach to the problem based on the conformal scalar
curvature theory was provided by R. Schoen \cite{Sc1}.\par We have to
notice that many mathematicians obtained various conditions for a
Riemannian manifold to be isometric to a sphere, ones used the
infinitesimal conformal transformations (see M. Obata \cite{Ob3}, C.C.
Hsiung-L.W Stern \cite{Hs-St}, K. Yano-T. Nagano \cite{Ya-Na}, ...)
and others gave conditions on the sectional curvatures and the Ricci
curvature or particular bounds for a certain eigenvalue of the
Laplacian on $(M,g)$ and the Ricci curvature (see S.Deshmukh and A. Al.Eid
\cite{De-Al}, ...).\par In this paper, we obtain a condition
involving the Ricci and the scalar curvature for $(M,g)$ to be
conformally equivalent to a sphere. This will be based on comparison
results for the Yamabe  invariant. The main ingredients are
symmetrization process and isoperimetric comparison results due to
P.Berard, G.Besson and S.Gallot \cite{B-B-G}, which have been
extended by the second author in \cite{Gam1} and \cite{Gam2}.

\section {Introduction and statement of main results}\label{1}
Let $(M,g)$ be a compact Riemannian manifold of dimension $n$
without boundary. We denote by $Ric_{g}$ its Ricci curvature,
$r_{0}$ the infinimum of $r(x)=\inf \{Ric_{g}(u,u),u\in
T_{x}M,\left| u\right|=1\}$, the least eigenvalue of $Ric_{g}$ on
the
tangent space $TM$, $R$ its scalar curvature and $d$ its diameter.\\
The isoperimetric profile of $(M,g)$ is defined by
\begin{equation}\label{iso1}
 h(s)=\inf \{\frac{vol\partial \Omega }{vol(M,g)},\,\, \Omega
\subset M \,\, \textit{ s.t }\,\,\frac{vol\Omega
}{vol(M,g)}=s\},\,\, s \in [ 0,1]
\end{equation}
where  $\Omega \subset M$ are smooth domains with regular
boundaries.\\ Let $Is(s)$ be the isoperimetric profile of the model
space: the unit sphere $S^{n}$ of $\mathbb{R}^{n+1}$ equipped with
its canonical metric, that is
\begin{equation}\label{iso2}
 Is(s)=\frac{ vol(\partial B(s))}{vol\text{
}S^{n}}
\end{equation}
where $B(s)$ is a geodesic ball of $S^{n}$ such that $\frac{vol B(s)}{vol S^{n}}= s.$\\
The following result is de to Bérard, Besson and Gallot \cite{B-B-G}:\\ If $(M,g)$ satisfies
\begin{equation}\label{1.3}
r_{0}d^{2}\geq (n-1)\varepsilon \alpha ^{2} \hspace{,2cm}\big(
\varepsilon \in \{-1,0,1\} \textit{  and } \alpha \in
\mathbb{R}_{+}\big),
\end{equation}then $\forall s \in [0,1]$
\begin{equation}\label{iso3}
 d \,h(s)\geq a(n,\varepsilon ,\alpha ) Is(s).\end{equation}
where $a(n,\varepsilon,\alpha )$ is a constant depending on $n$,
$\varepsilon $ and $\alpha $ as follows
\begin{equation}\label{a}
 a(n,\varepsilon ,\alpha
)=\left\{\begin{array}{lll} \alpha \sigma
_{n}^{\frac{1}{n}}\big[2\int_{0}^{\frac{\alpha }{2}}(\cos
t)^{n-1}dt\big]^{-\frac{1}{n}},\textit{ if }\varepsilon =1 &  &  \\
(1 + n\,\sigma _{n})^{\frac{1}{n}} - 1,\textit{ if }\varepsilon =0 &  &  \\
\alpha c(\alpha ),\textit{ if }\varepsilon = -1 ,
\end{array}\right.
\end{equation}
 where $\sigma _{n}=\int_{0}^{\pi }(\sin t)^{n-1}dt$ and $c(\alpha )$ is the unique root of the equation $\sigma _{n}(\cosh
y)^{n}=\sinh y\int_{y}^{y+\alpha }(\cosh t)^{n-1}dt$.
 The Inequality  (\ref{1.3}) is sharper than the one given by M. Gromov \cite{Gr}.\\
Let $(M,g)$ be a  compact Riemannian manifold without boundary
satisfying relation (\ref{1.3}). In \cite{Gam2} ( Theorem5
), the first author obtained a lower bound for the least eigenvalue of the operator $\Delta
+C$, where $\Delta $ is the Laplacian of $(M,g)$ and $C$ is a
potential.\\\par Following the proof of \cite{Gam2}, we give in the
first step a lower bound for the least eigenvalue  of the conformal
Laplacian of $(M,g)$, $L = c_n\Delta + R $  with $c_n=
4\frac{n-1}{n-2}$.

 We begin by providing some notations: let $V$ denote the
volume of $(M,g)$ and $\omega _{n}$ the one of $(S^{n},can)$. For
given positive reals $r_{1},r_{2}$ $(0<r_{1},r_{2}<\pi )$, let
$B(S,r_{1})$, ( resp. $B(N,r_{2})$ ) be the geodesic ball of $S^{n}$
of center the south pole $S$ and radius $r_{1}$ (resp. the geodesic
ball of $S^{n}$ of center the north pole $N$ and radius $r_{2})$.

\begin{proposition}\label{oo} Let  $(M,g)$ be a compact Riemannian
manifold of dimension $n$ without boundary satisfying relation
(\ref{1.3}). Let $\mu_{1}(M)$ (resp $\rho _{1}(S^{n})$) denotes the
least eigenvalue of the conformal laplacian $L = c_{n}\Delta + R $
on $\mathit{M}$ (resp. of $c_{n}\Delta + h_{+} - h_{-} $ on the unit
sphere $S^{n})$ acting on functions with
\begin{equation}\label{h+}
\mathit{h}_{+}\mathit{=}\left\{
\begin{array}{lll}
\text{ }(d/a)^{2} \sup R_{+} \qquad \text{on }B(S,r_{1}) \\\\
0 \qquad  \qquad \qquad \qquad  \text{on the complementary in }S^{n}
\end{array}
\right.\end{equation}
 and\begin{equation}\label{h-}
\mathit{h}_{-} = \left\{
\begin{array}{lll}
\text{ }(d/a)^{2} \sup R_{-} \qquad\text{on } B(N,r_{2}) \\\\
0 \qquad  \qquad \qquad \qquad \text{on the complementary in }
S^{n},
\end{array}
\right.
\end{equation}
where $r_{1},\; r_{2}$ satisfy:
\begin{equation}\label{r1}
\omega _{n}^{-1}vol B(S,r_{1} ) = V^{-1} \| R_{+} \|_{L^{1}(M)}\,/
\,\| R_{+} \| _{L^{\infty }(M)}
\end{equation}
and
\begin{equation}\label{r2}
\omega _{n}^{-1} vol B(N,r_{2}) = {V}^{-1}\| R_{-} \|_{L^{1}(M)}\,/
\,\| R_{-}\|_{L^{\infty }(M)}.
\end{equation}
Then
\begin{equation}\label{prop1}
\mu_{1}(M) \geq (a/d)^{2}{\rho }_{1}(S^{n}).
\end{equation}
\end{proposition}
\vspace*{0,5cm} In the second step, we give a lower bound for the
first Yamabe invariant of $(M,g)$ which we denote by $\lambda(M)$.
Let $\rho(S^{n})$ be the least eigenvalue on $S^{n}$ of
\begin{equation}
 c_n\Delta u + h u = \rho
u^{\frac{n+2}{n-2}}
\end{equation}
where $h=h_+ -h_-$ is given by (\ref{h+}) and (\ref{h-}), we have
\begin{theorem}\label{o1} Let $(M,g)$ be a compact Riemannian
manifold of dimension $n$ without boundary satisfying relation
(\ref{1.3}). We have
\begin{equation}
\lambda(M)\geq (a/d)^{2}\beta^{\frac{2}{n}}\rho(S^{n})
\end{equation}
\end{theorem}
Since the Yamabe invariant of  $(M,g)$ is bounded from above by the
one of the sphere, we derive the following rigidity result.
\begin{theorem}\label{o2} A compact Riemannian manifold $(M,g)$ of
dimension $n\geq 3$ without boundary satisfying condition (\ref{1.3})
with $\varepsilon =1$, and \begin{equation}\label{ine1}
(a/d)^{2}\beta^{\frac{2}{n}}\rho(S^{n})\geq\lambda(S^{n})
\end{equation} is conformally diffeomorphic to the unit sphere $S^{n}$
\end{theorem}
Notice that in the particular case where the scalar curvature is
constant, we have
$\lambda(M)=RV^{\frac{2}{n}}=(\frac{a}{d})^{2}\beta^{\frac{2}{n}}\rho(S^{n})$
and condition (\ref{ine1}) is only satisfied when M is conformally
diffeomorphic to the sphere.\\
 As a consequence of Proposition \ref{oo} and Theorem \ref{o1} and in the case where $Ric\geq n-1$, we obtain the following
result:
\begin{corollary}\label{o3} Let $(M,g)$ be a compact Riemannian manifold of
dimension $n\geq 3$ without boundary satisfying $\ Ric\geq n-1$.
Then
\begin{equation}\label{coro1}
\mu_{1}(M)\geq n(n-1) = \mu_{1}(S^{n})
\end{equation}
\begin{equation}\label{coro2}
\lambda(M)\geq n(n-1)V^{\frac{2}{n}} = \big(
\frac{V}{\omega_{n}}\big)^{\frac{2}{n}}\lambda(S^{n}).
\end{equation}
Where $\lambda(M)$ denote the Yamabe invariant of $(M,g)$ and
$\lambda(S^{n})$ the one of $S^{n}$.\end{corollary}
 We have to point out that (\ref{coro1}) and (\ref{coro2}) are optimal in the
case where the metric $\mathit{g}$ is Einstein.\\
The inequality (\ref{coro2}) was proved by J.Petean and S.Ilias in
\cite{Pe} and \cite{Il} respectively by using analogous methods.
\section {Yamabe Problem and conformal invariant $\lambda(M)$ }\label{2}
 {\bf Yamabe question:} Given a compact Riemannian
manifold $(M,g)$ without boundary of dimension $n\geq 3,$ is there a
metric $\widetilde{g}$ conformal to $g$ which has constant scalar
curvature $R_{\widetilde{g}}=\lambda $?. We write
$\widetilde{g}=u^{\frac{4}{n-2}}g$, $u>0$. By a simple computation
we obtain:
\begin{equation}\label{yr}
 R_{\widetilde{g}}=u^{- \frac{n+2}{n-2}}(
c_{n} \Delta u + R u)
\end{equation}
where $R$ is the scalar curvature and $\Delta $ the Laplacian of
$(M,g)$. Hence the Yamabe problem is equivalent to solve:
\begin{equation}\label{ye1}
 c_{n}\Delta u + R u = \lambda
u^{\frac{n+2}{n-2}},\,\, u > 0
\end{equation}
We will use the following notations: $p=\frac{2n}{n-2}$ and $L=c_{n}\Delta + R$.\\
The operator $L $ is called the conformal laplacian of $(M,g)$.
Equation (\ref{ye1}) can be rewritten as
\begin{equation}\label{ye2}
 Lu=\lambda u^{p-1},\,\, u > 0.
\end{equation}
Yamabe noticed that (\ref{ye1}) is the Euler-Lagrange equation of
the functional:
\begin{equation}\label{fun}
Q_{0}(\widetilde{g})=\frac{\int_{M}R_{\widetilde{g}} d
v_{\tilde{g}}}{\big( \int_{M}dv_{\tilde{g}}\big) ^{\frac{2}{p}}}
\end{equation}
when restricted to a conformal class $[g]=\{hg \,/ \, h\in C^{\infty
}(M),\, h > 0\}$, where $dv_{\tilde{g}}$ is the volume form of
$(M,\widetilde{g})$ and $h = u^{p-2}$, $u > 0$. In fact, on $[g]$ we
can write $Q_{0}(\widetilde{g}) = Q_{0}(u^{p-2}g) = J(u)$, where
\begin{equation}\label{j}
 J(u) = \frac{\int_{M} u L u d
v_{g}}{\big(\int_{M}u^{p}dv_{g}\big)^{\frac{2}{p}}} =
\frac{\int_{M}(c_{n}|\nabla u|^{2} + R u^{2})d v_{g}}{\| u
\|_{p}^{2}}.
\end{equation}
We call $J(u)$ the Yamabe quotient of $(M,g)$. Let $u$ be a positive
function in $ C^{\infty }(M)$ and  a critical point of $J$, then it
is easy to see that $u$ satisfies equation (\ref{ye1}) with
$\lambda = J(u)$.\\
By using a H\"{o}lder inequality, we derive that the functional $J$
is bounded from below. The infimum
\begin{eqnarray}\label{yinv1}
\lambda (M)= \inf \big\{Q_{0}(\widetilde{g})\,/\,\widetilde{g}\in
[g]\big\} = \inf \big\{J(u)\,/\, u\in C^{\infty }(M)\text{, }u >
0\big\}
\end{eqnarray}
is a conformal invariant, which means that it is
determined by the conformal class and is independent of the choice
of the initial metric $g$ in the conformal class.
It is called the Yamabe invariant of $(M,g)$. \\
We have the following results
\\\textbf{Theorem A:}(\cite{Ya},
\cite{Tr},\cite{Au1},\cite{Au2}):  \emph{The Yamabe problem can be
solved on any compact manifold $M$  with }$\lambda (M) < \lambda
(S^{n})$.\\
\textbf{Theorem B:} (\cite{Ya}, \cite{Au1}, \cite{Au2}):\emph{ For
any compact Riemannian manifold $(M,g)$ without boundary, we always
have} $\lambda (M) \leq \lambda (S^{n}) = n(n-1)\omega
_{n}^{\frac{2}{n}}$.\vspace*{0,2cm}\par Theorem A reduces the
resolution of Yamabe problem to the estimate of the invariant
$\lambda (M)$. In fact, if we can find a function $u \in
L_{1}^{2}(M)$ such that $J(u) < \lambda (S^{n})$, then $\lambda (M)
< \lambda (S^{n})$, hence the Yamabe problem has a solution. \\In
this way T.Aubin \cite{Au2} proved the conjecture in the two
following cases:\\1) $(M,g)$ is not a conformally flat compact
Riemannian manifold of dimension $n\geq 6$.\\2) $(M,g)$ is a locally
conformally flat compact Riemannian manifold of dimension $n\geq 3$
and finite Poincar\'{e} group, not conformal to $(S^{n},can)$.\\
R. Schoen \cite{Sc2} solved all the remaining cases of the Yamabe
problem,
using the positive mass theorem. \\
We remark that for the case where $(M,g)$ is conformal to $S^{n}$,
the Yamabe problem clearly has a solution. If $\Phi :\, M
\rightarrow  S^{n}$ is a conformal diffeomorphism then
$\Phi^{*}(g_{0}) = f\,g$, where $g_{0}$ is the standard metric of
$S^{n}$ and $f$ a positive function in  $C^{\infty }(M)$, clearly
$fg$ has constant scalar curvature.\\ Besides the proof of T.Aubin
and R.Schoen of the Yamabe problem, another proof by A.Bahri
\cite{Ba}, A.Bahri-H.Br\'{e}zis \cite{Ba-Br} of the same conjecture
is available using the theory of critical points
at infinity.
\section { Symmetrization method and applications }\label{3}

In the following
we give lower bounds for the first eigenvalue of the conformal
Laplacian of the manifold $(M,g)$, that we denote by $\mu_{1}(M)$
and its Yamabe invariant $\lambda (M)$. The method we use here is
inspired by the one used
in \cite{Gam1} and \cite{Gam2}.\\
We begin by the case of the least eigenvalue of the Laplacian and the proof
obtain Proposition \ref{oo}.
\begin{proof}[Proof of Proposition \ref{oo}]: Let $(M,g)$ be a compact Riemannian manifold which satisfies the
  isoperimetric inequality (\ref{1.3}). One can apply the symmetrization process described
  in \cite{Be1}, \cite{Be2} and \cite{Ta} or \cite{Gam1} and \cite{Gam2}
  to symmetrize a smooth function $f$ into a radial function $f^{*}$ on the model space $(S^{n},can)$.
   The function $f^{*}$  is in $H^{1}(S^{n})$, radial ( w.r.t the north pole ) and
   satisfies
 \begin{eqnarray}\label{coaire}
\left\{
\begin{array}{lll}\omega _{n}\int_{M}f^{q}dv_{g}&=&V\int_{S^{n}}f^{*}\; ^{q}dv,\text{
for all real } q \geq 1
\\\\\omega _{n}\int_{M} |\nabla f |^{2}dv_{g} &\geq &
V(\frac{a}{d})^{2} \int_{S^{n}} \mid \nabla f^{*} \mid^{2}dv
\end{array} \right.
\end{eqnarray}
The first identity of (\ref{coaire}) derives from the coarea
formula (see \cite{Ban})  and the second
can be proved through coarea formula, isoperimetric inequality of \cite{B-B-G} and the Cauchy-Schwarz inequality .\\
The inequality of Hardy-Littlewood-Polya ( \cite{Ta} formulas $(60)$
and $(13))$ implies
\begin{equation}\label{hardy}
\int_{M}R f^{2}dv_{g}\geq \beta
\int_{0}^{\omega_{n}}\big[R_{+}^{*}(V-\beta u) - R_{-}^{*}(\beta
u)\big]f^{* 2}(u)du
\end{equation}
where $\beta $ denotes the ratio $\frac{V}{\omega _{n}}$,
$R_{+}^{*}(V - \beta u)$ is the the increasing symmetric
rearrangement of $R_{+}$ and $R_{-}^{*}(\beta u)$ (respectively
$f^{* }(u)$) the decreasing symmetric rearrangement of\;
 $R_{-}$\;(respectively of $f$).
 Then we apply at the right handside of (\ref{hardy})
 the following Steffensen inequality (one can see D.S Mitrinovi\'{c}
 \cite{Mi}):

 \bigskip

\textbf{Theorem} \label{M}(\cite{Mi}) \emph{Let} $\varphi $
\emph{and} $\psi $ \emph{be two given integrable functions defined
on the interval} $(a,b)$ \emph{such that} $\varphi $ \emph{is
decreasing and} $0\leq \psi \leq 1$ \emph{on} $(a,b)$, \emph{then}:
$$
\int_{b-\gamma }^{b}\varphi (t)dt\leq \int_{a}^{b}\varphi
(t)\psi(t)dt \leq \int_{a}^{a+\gamma }\varphi (t) dt,$$
\emph{where}
 $\gamma =\int_{a}^{b} \psi (t) dt$.\\
 We obtain
\begin{equation*}
\int_{M}R f^{2}dv_{g}\geq \beta \Big[\sup R_{+}^{*}\int_{\gamma
_{+}}^{\omega _{n}}f^{* 2}(u)du-\sup R_{-}^{* }\int_{0}^{\gamma
_{-}}f^{* 2}(u)du\Big]
\end{equation*}
with
\begin{eqnarray*}\gamma _{+} &=& \omega _{n} - (\beta \sup
R_{+}^{*})^{-1}\int_{0}^{V}R_{+}^{*}(V - u)d u \\and\\\gamma _{-}
&=& (\beta \sup R_{-}^{*})^{-1}\int_{0}^{V}R_{-}^{*}(u) d u
\end{eqnarray*}

We identify $R_{+}^{*}$ (respectively  $R_{-}^{*}$ ) with a function
 $R_{+}^{*}(\pi-r)$ (respectively $R_{-}^{*}(r)$) of the distance to the north pole.
Let $\widetilde{R}_{+}$ and $\widetilde{R}_{-}$  be the radial functions
defined on $S^{n}$ as follows
\begin{equation}\label{R+}
\qquad\qquad\quad \widetilde{R}_{+}(r) = \left\{
\begin{array}{lll}
(d/a)^{2}  \sup R^{*}_{+}  \quad on \quad [\pi-r_{1},\pi] \\\\
0 \qquad  \qquad  \qquad on\;\; [0,\pi-r_{1}]
\end{array}\right.
\end{equation}and
\begin{equation}\label{R-}
\qquad\qquad\quad \widetilde{R}_{-}(r) = \left\{
\begin{array}{lll}
(d/a)^{2}  \sup R^{*}_{-}  \quad\quad on\quad [0 , r_{2}] \\\\
0 \qquad  \qquad \qquad  \qquad on \;\; [r_{2} , \pi].
\end{array}\right.
\end{equation}
 We have
$$(d/a)^{2}\sup R_{+}^{*}\int_{\gamma
_{+}}^{\omega _{n}}f^{* 2}(u)du=
\omega_{n-1}\int_{\pi-r_{1}}^{\pi}\widetilde{R}_{+}(r)
f^{*2}(r)(\sin r)^{n-1}dr$$and
$$(d/a)^{2}\sup R_{-}^{* }\int_{0}^{\gamma
_{-}}f^{* 2}(u)du=\omega_{n-1}\int_{0}^{r_{2}}\widetilde{R}_{-}(r)
f^{*2}(r)(\sin r)^{n-1}dr.$$ Therefore, we obtain
\begin{equation}
  \int_{M}R f^{2} d v_{g} \geq \beta
(\frac{a}{d})^{2}\int_{S^{n}}(h_{+} - h_{-})f^{* 2}(v)dv
\end{equation}
and finally using (\ref{coaire}),
\begin{equation}
  \frac{\int_{M}f L f d v_{g}}{\int_{M}f^{2}dv_{g}}\geq
(\frac{a}{d})^{2}\frac{\int_{S^{n}}\big[c_{n}| df^{*}|^{2} + (h_{+}
- h_{-}) f^{* 2}\big] dv}{\int_{S^{n}}f^{* 2}dv},
\end{equation}
Hence we end the proof by using the fact that the least eigenvalue
is the infimum of the Rayleigh quotient.\end{proof}
\par In the sequel we deal with the Yamabe invariant $\lambda
(M)$ of $(M,g)$ introduced in (\ref{yinv1}). Since this invariant
can be expressed in terms of Rayleigh quotient as
\begin{equation}
\lambda (M) = \inf\limits_{u}\frac{\int_{_{M}}\big( c_{n}|\nabla
u|^{2} + R u^{2}\big) dv_{g}}{\big( \int_{M}u^{p}
dv_{g}\big)^{\frac{2}{p}}}
\end{equation}
\noindent where the infinimum is taken over all smooth real-valued
 positive functions $u$ on $M$, we can use the same techniques
introduced above in the aim to give bounds for $\lambda (M)$. Let
$\rho(S^{n})$ be the least eigenvalue on $S^{n}$ of
\begin{equation}\label{ro}
 c_n\Delta u + h u = \rho
u^{\frac{n+2}{n-2}}\qquad\qquad\qquad\qquad\qquad
\end{equation}
where $h=h_+ -h_-$ is given by (\ref{h+}) and (\ref{h-})


\begin{proof}[Proof of Theorem \ref{o1}]: We begin the proof by providing a lower
bound for the Yamabe invariant $\lambda (M)$ with the use of the
symmetrization method given in Proposition \ref{oo}. For a positive
function $f$ in $C^{\infty }(M)$, we consider its decreasing
symmetric rearrangement $f^{*}$. Let $R$ be the scalar curvature of
$(M,g)$ and $h$ the function defined on the unit sphere by
(\ref{h+}) and (\ref{h-}).
Following the same steps as in the proof of Proposition \ref{oo},
and applying (\ref{coaire}) for $q=p$, we obtain
\begin{equation}
\frac{\int_{_{M}}\big( c_{n}|\nabla f|^{2} + R f^{2}\big)
dv_{g}}{\big( \int_{M} f^{p} dv_{g}\big)^{\frac{2}{p}}}\geq
\frac{\beta }{\beta^{\frac{2}{p}}} \big(
\frac{a}{d}\big)^{2}\frac{\int_{_{S^{n}}} \big( c_{n}|\nabla
f^{*}|^{2} +  h f^{^{*}2}\big) dv}{\big( \int_{_{S^{n}}}f^{* p} dv
\big)^{\frac{2}{p}}},
\end{equation}
and hence\begin{equation}\label{ine} \lambda(M)\geq
(a/d)^{2}\beta^{\frac{2}{n}}\rho(S^{n}).
\end{equation}\end{proof}

\begin{proof}[Proof of Theorem \ref{o2}] From Theorem \ref{o1},
 we have
\begin{eqnarray*}
\lambda(M)\geq (a/d)^{2}\beta^{\frac{2}{n}}\rho(S^{n}).
\end{eqnarray*}
On the other hand  $\lambda(M)$ is
upper bounded by $\lambda(S^{n})$, and since in this case $R_{-}=0$,
we derive that
$(a/d)^{2}\beta^{\frac{2}{n}}\rho(S^{n})\leq\lambda(S^{n})$. Hence,
if $(a/d)^{2}\beta^{\frac{2}{n}}\rho(S^{n})\geq\lambda(S^{n})$, then the equality $\lambda(M)=\lambda(S^{n})$ holds and
$(M,g)$ is conformally  diffeomorphic to the unit sphere $S^{n}$,
which completes the proof.\end{proof}
 As a consequence of Proposition \ref{oo} and Theorem \ref{o1} in the
 case where
$Ric\geq n-1$ ( take $\alpha = d$ and $\varepsilon = 1$ in
(\ref{1.3})), we obtain the corollary given in the introduction.
\begin{proof}[Proof of corollary \ref{o3}]  Since the scalar curvature $R\geq
n(n-1)$ we derive the following inequalities
$$\inf_{f}\frac{\int (c_{n}|\nabla f|^{2}+ R f^{2})dv_{g}}{\int f^{2}dv_{g}}\geq
 \inf_{f}\frac{\int( c_{n}|\nabla f|^{2}+n(n-1)f^{2})dv}{\int f^{2}dv}$$
 and  $$\inf_{f}\frac{(\int c_{n}|\nabla f|^{2}+ R f^{2})dv_{g}}{(\int f^{p}dv_{g})^{\frac{2}{p}}}\geq
 \inf_{f}\frac{\int( c_{n}|\nabla f|^{2}+n(n-1) f^{2})dv}{(\int f^{p}dv)^{\frac{2}{p}}},$$
 where the infinimum is taken over all smooth real-valued positive functions $f$
 on $M$. The same arguments as in the proof of Proposition \ref{oo} and
Theorem \ref{o1} enable us to lowerbound the right hand sides of
these inequalities by $\mu_{1}(S^{n})$ and $\beta ^{\frac{2}{n}}
\lambda (S^{n})$ respectively.\end{proof}


\end{document}